\documentclass[12pt,a4paper]{article}
\usepackage[latin2]{inputenc}
\usepackage{amsmath}
\usepackage{amsfonts}
\usepackage{amssymb}
\usepackage{graphicx}
\usepackage[left=2.00cm, right=2.00cm, top=2.50cm, bottom=2.50cm]{geometry}
\usepackage[T1]{fontenc}

\def\Q{{\mathbb Q}}
\def\Z{{\mathbb Z}}

\def\C{{\mathbb C}}

\newtheorem{rem}{Remark}
\newtheorem{thm}{Theorem}[section]

\newtheorem{lem}[thm]{Lemma}

\def\ve{\varepsilon}

\title{
Monogenity of pure quintic fields: the power of sieving
}
\author{
Istv\'an Ga\'al\\
{\small University of Debrecen, Mathematical Institute} \\
{\small H--4002 Debrecen Pf.400., Hungary,} \\
{\small e--mail: gaal.istvan@unideb.hu},
}

\begin{document}
\baselineskip=17pt

\maketitle
\thispagestyle{empty}

\renewcommand{\thefootnote}{\arabic{footnote}}
\setcounter{footnote}{0}

\vspace{0.5cm}


\begin{abstract}
We provide a simple algorithm for calculating all generators of power integral bases
in pure quintic fields. This procedure involves the usual standard elements like
Baker's method, LLL-reduction. The main purpose of the paper is to 
introduce a new idea to considerably diminish the number of small exponents to be considered
after the reduction step.
This new idea allows to test all remaining small exponents within a few minutes, 
using an appropriate sieve method, which turns out to be surprisingly fast. 
This idea will be applicable in many similar cases.
\end{abstract}

\section{Introduction}

Monogenity and power integral bases is a classical topic in algebraic number theory,
which is intensively studied even today, see \cite{nark} for classical results,
\cite{book} and \cite{axioms} for more recent results.

A number field $K$ of degree $n$ with ring of integers $\Z_K$ 
is called {\it monogenic} (cf. \cite{book}) if there 
exists $\xi\in \Z_K$ such that $(1,\xi,\ldots,\xi^{n-1})$ is an integral basis, 
called {\it power integral basis}. 
We call $\xi$ the {\it generator} of this power integral basis.

An irreducible polynomial $f(x)\in\Z[x]$ is called {\it monogenic}, 
if a root $\xi$ of $f(x)$ generates a power integral basis in $K=\Q(\xi)$.
If $f(x)$ is monogenic, then $K$ is also monogenic, but the converse is
not true.

For $\alpha\in\Z_K$ (generating $K$ over $\Q$)  the module index 
\[
I(\alpha)=(\Z_K:\Z[\alpha])
\]
is called the {\it index} of $\alpha$. The element 
$\alpha$ generates a power integral basis in $K$
if and only if $I(\alpha)=1$. 
The elements $\alpha,\beta\in\Z_K$ are called {\it equivalent}, if
$\alpha\pm \beta\in\Z$. Equivalent algebraic integers have equal indices.
Searching for elements of $\Z_K$, generating power integral bases, leads to
a Diophantine equation, called {\it index form equation} (cf. \cite{book}). 

There are efficient algorithms for determining generators of power integral bases
in cubic and quartic number fields (cf. \cite{book}), because in those cases the corresponding
index form equation can be reduced to Thue equations. The algorithms for solving absolute
(over $\Z$) and relative (over $\Z_K$) Thue equations were also applied to solve
index form equations in special higher degree fields, like sextic and octic fields with
a quadratic subfield (cf. \cite{book}). 

There were attempts to solve index form
equations in general quintic \cite{quintic} and sextic \cite{sextic} 
number fields, but they turned out to be unusable, because of a very long CPU time
(8 hours in the quintic case, 2 months in the sextic case, which is very long, 
even if we encounter the improvement of hardware since that time).
This showed that the method used to derive effective results for index form equations
can not be used in practice, not even if it is together with an efficient enumeration
algorithm (the so called ellipsoid method), developed by Wildanger \cite{wild}
and Ga\'al and Pohst \cite{gp}. 

On the other hand, using Newton polygon method and Dedekind's criterion 
many papers appeared about the monogenity of several types of pure fields
and binomials $x^n-m$ (cf. \cite{axioms}). Parallelly, we started to implement 
the known methods to solve index form equations in certain pure fields.
This was successful in pure sextic (cf. \cite{p6}, \cite{p62}, \cite{sieve}) 
and pure octic fields (cf. \cite{p8}). The reason is that in those cases the
pure field has a quadratic subfield and the index form has a factor which implies a relative
Thue equation.

It has also turned out (cf. \cite{sieve}), that the sieve method (see \cite{book})
can be surprisingly efficient in such calculations.

In the present paper we develop an algorithm for calculating all inequivalent generators
of power integral bases in pure quintic fields. Our method involves the usual standard 
ingrediendts like Baker's method, LLL-reduction. The most difficult and time consuming 
part in \cite{quintic} was the enumeration of small exponent tuples, remaining after
the reduction procedure. Here we replace the complicated ellipsoid method used in \cite{quintic}
by a simple sieve method. However, in the case of pure quintic fields, the sieve method 
is only sufficient if we use an idea, applied first time in this paper.
This is the main novelty of the paper. 
This new idea allows us to present all generators of power integral bases in 
pure quintic fields within about 5 minutes of CPU time.
It is certain that this idea will have many further applications.

\section{Pure quintic fields}

Let $m\in\Z$, such that $f(x)=x^5-m$ is irreducible.
For simplicity we also assume that it is also monogenic.
For a straightforward extension of the method to non-monogenic polynomials see Section \ref{conc}.

The purpose of the present paper is to find all possible generators of power integral bases
in $K=\Q(\xi)$, where $\xi$ is a root of $f(x)$. 
We  wonder if there are any
further inequivalent generators of power integral bases in $K$, in addition to $\xi$.

If $f(x)$ is monogenic, we can represent any $\alpha\in\Z_K$ in the form
\begin{equation}
    \alpha=x_0+x_1\xi+x_2\xi^2+x_3\xi^3+x_4\xi^4,
\label{ax}
\end{equation}
with $x_0,x_1,x_2,x_3,x_4\in\Z$. Denote by $D_K$ the discriminant of $K$, 
and let $\gamma^{(j)}$ be the conjugates of any $\gamma\in K$ ($j=1,\ldots,5$).
We have
\[
\prod_{1\le i< j\le 5}|\alpha^{(i)}-\alpha^{(j)}|=I(\alpha)\cdot \sqrt{|D_K|},
\]
and
\[
\prod_{1\le i< j\le 5}|\xi^{(i)}-\xi^{(j)}|=I(\xi)\cdot \sqrt{|D_K|},
\]
whence
\begin{equation}
  \prod_{1\le i< j\le 5}\frac{|\alpha^{(i)}-\alpha^{(j)}|}{|\xi^{(i)}-\xi^{(j)}|}  =
  \frac{I(\alpha)}{I(\xi)}.
\label{alpha}
\end{equation}
Using the representation (\ref{ax}) of $\alpha$ this can be written as
\begin{equation}
\prod_{1\le i< j\le 5}|\ell_{ij}(x)|=\frac{I(\alpha)}{I(\xi)},
\label{xij}
\end{equation}
where
\[
\ell_{i,j}(x)=x_1+x_2(\xi^{(i)}+\xi^{(j)})
+x_3((\xi^{(i)}+\xi^{(j)})^2-\xi^{(i)}\xi^{(j)})
\]
\[
+x_4((\xi^{(i)}+\xi^{(j)})^3-2\xi^{(i)}\xi^{(j)}(\xi^{(i)}+\xi^{(j)})).
\]
The Galois group of $f(x)$ is $D_5=ST3$ which is doubly transitive (cf. \cite{quintic}).
Denote by $L_{i,j}=Q(\xi^{(i)}+\xi^{(j)},\xi^{(i)}\xi^{(j)})$ the subfield
of $Q(\xi^{(i)},\xi^{(j)})$ which remains fixed under the action 
$(i,j)\to (j,i)$ of the Galois group. This is a proper subfield of $Q(\xi^{(i)},\xi^{(j)})$
of degree 10. We denote by $\gamma^{(i,j)}$ the conjugates of any $\gamma^{(1,2)}\in L_{1,2}$
corresponding to $\xi^{(i)}+\xi^{(j)},\xi^{(i)}\xi^{(j)}$. Obviously,
$\gamma^{(i,j)}=\gamma^{(j,i)}\; (1\le i<j\le 5)$.

Equation (\ref{xij}) implies that $\ell_{i,j}(x)$ is contained in $L_{i,j}$.
In our case $I(\xi)=1$ and we are looking for $\alpha\in\Z_K$ with $I(\alpha)=1$, therefore
$\ell_{i,j}$ is a unit in $L_{i,j}$. The field $L_{i,j}$ has $s=2$ real and $2t=8$ 
complex conjugates, therefore its unit rank is $r=s+t-1=5$. 
Let $\ve_1^{(i,j)},\ldots,\ve_5^{(i,j)}$ be a set of fundamental units in $L_{i,j}$.
We have 
\begin{equation}
\ell_{i,j}(x)=\pm \prod_{h=1}^5 (\ve_h^{(i,j)})^{a_h}
\label{spd}
\end{equation}
with $a_h\in\Z\; (1\le h\le 5)$.

\section{The structure of the degree 10 field $L_{i,j}$}
\label{str}

Using symmetric polynomials it is easily seen that
the field
$L_{i,j}=Q(\xi^{(i)}+\xi^{(j)},\xi^{(i)}\xi^{(j)})$
is generated by a root of
\[
g(x)=x^{10}+11mx^5-m^2.
\]
Denote by $\beta$ a root of $g(x)$. Then 
\[
(\beta^5)^2+11m\beta^5-m^2=0,
\]
whence
\[
\beta^5=\frac{-11m\pm\sqrt{121m^2+4m^2}}{2}=m\frac{-11\pm5\sqrt{5}}{2},
\]
which implies that $\sqrt{5}\in L_{i,j}$.
Moreover,
\[
\beta^5=(\sqrt[5]{m})^5\cdot (-1)^5\cdot \left(\frac{1\mp \sqrt{5}}{2}\right)^5,
\]
whence a value of $\sqrt[5]{m}$ is also contained in $L_{i,j}$.
These are the quadratic and quintic subfields of $L_{i,j}$.
Let $\eta$ be the fifth root of unity with smallest positive argument,
that is $\eta\approx \sqrt[5]{1}=-0.8090-0.5877 I$ and let 
$\sqrt[5]{m}$ be the real value of the fifth root of $m$.
Then $\xi^{(i)}=\eta^i\cdot \sqrt[5]{m}\; (1\le i\le 5)$.
We observe that 
\[
\frac{(\xi^{(i)}+\xi^{(f)})^2}{\xi^{(i)}\xi^{(j)}}
=\frac{\eta^{2i}+\eta^{2j}+2\eta^{i+j}}{\eta^{i+j}}
=\eta^{i-j}\eta^{j-i}+2.
\]
This value is equal to 
\[
f_1=\frac{3+\sqrt{5}}{2}\;\; \text{for}\;\; (i,j)=(1,3),(1,4),(2,4),(2,5),(3,5)
\]
and it is
\begin{equation}
f_2=\frac{3-\sqrt{5}}{2}\;\; \text{for}\;\; (i,j)=(1,2),(1,5),(2,3),(3,4),(4,5).
\label{f2}
\end{equation}
This way we can detect in which $L_{i,j}$ the elements 
of the quadratic subfield $\Q(\sqrt{5})$ remain fixed.
This will be important in the final step of our algorithm.

\section{Baker's method}

The first ingredients of our algorithm are the standard ones, Baker's method
and reduction. We following \cite{quintic}, but make that 
description more exact in our case.

Let $1\le i,j,k\le 5$ be distinct indices. We have
\begin{equation}
(\alpha^{(i)}-\alpha^{(j)})+(\alpha^{(j)}-\alpha^{(k)})+(\alpha^{(k)}-\alpha^{(i)})=0
\end{equation}
whence
\begin{equation}
(\xi^{(i)}-\xi^{(j)})\ell_{i,j}(x)+(\xi^{(j)}-\xi^{(k)})\ell_{j,k}(x)
+(\xi^{(k)}-\xi^{(i)})\ell_{k,i}(x)=0.
\end{equation}
that is 
\begin{equation}
\frac{(\xi^{(i)}-\xi^{(j)})\ell_{i,j}(x)}{(\xi^{(i)}-\xi^{(k)})\ell_{i,k}(x)}   
+
\frac{(\xi^{(k)}-\xi^{(j)})\ell_{j.k}(x)}{(\xi^{(k)}-\xi^{(i)})\ell_{i,k}(x)} 
=1.
\label{xxijk}
\end{equation}
This can be written as
\begin{equation}
\delta_{ijk}\mu_{ijk}  + \delta_{kji}\mu_{kji}  =  1,
\end{equation}
where
\[
\delta_{ijk}=\frac{\xi^{(i)}-\xi^{(j)}}{\xi^{(i)}-\xi^{(k)}},\;\; 
\mu_{ijk}=\pm \prod_{h=1}^5 \left(\frac{\ve_h^{(i,j)}}{\ve_h^{(i,k)}}\right)^{a_h}
\]
We have
\[
a_1\log\left|\frac{\ve_1^{(i,j)}}{\ve_1^{(j,k)}}\right|
+\ldots+
a_5\log\left|\frac{\ve_5^{(i,j)}}{\ve_5^{(j,k)}}\right|
=
\log\left|\mu_{ijk}\right|.
\]
Select five distinct triples $(i,j,k)$ for which the coefficient matrix $T$ of this system
of linear equations is regular. Denote by $c_1$ the row norm 
(the maximum sum of the absolute values of entries in a row) of $T^{-1}$.
Then we obtain 
\[
A=\max_{1\le h\le 5}|a_h|\le c_1\left|\log\left|\mu_{i_0j_0k_0}\right|\right|
\]
where
\[
\left|\log\left|\mu_{i_0j_0k_0}\right|\right|
=\max_{i,j,k}\left|\log\left|\mu_{ijk}\right|\right|.
\]
If $\log\left|\mu_{i_0j_0k_0}\right|<1$, then this implies
\[
\left|\mu_{i_0j_0k_0}\right|\le \exp \left(-\frac{A}{c_1}\right),
\]
otherwise the same holds for $\mu_{i_0k_0j_0}=1/\mu_{i_0j_0k_0}$.

We omit the subindices of $i_0,j_0,k_0$ and assume that 
\[
\left|\mu_{ijk}\right|\le \exp \left(-\frac{A}{c_1}\right).
\]
\begin{rem}
As we cannot predict, for which triple $i,j,k$ this inequality is satisfied,
we have to apply the following estimates and the reduction procedure for all
possible triples $i,j,k$. After having the reduced bound for $A$ in all possible
cases, the enumeration procedure must be performed only once.
In the enumeration we have to use the maximum of the reduced bounds obtained
for the possible triples $i,j,k$.
\end{rem}

We have
\[
\Lambda=\left|
\log\left|\frac{\xi^{(k)}-\xi^{(j)}}{\xi^{(k)}-\xi^{(i)}}\right|
+a_1\log\left|\frac{\ve_1^{(j,k)}}{\ve_1^{(i,k)}}\right|+\ldots 
+a_5\log\left|\frac{\ve_5^{(j,k)}}{\ve_5^{(i,k)}}\right|
\right|
\]
\[
=\left|\log\left|\delta_{kji}\mu_{kji}\right|\right|
\le 2 \left|1-\left|\delta_{kji}\mu_{kji}\right|\right|
\le 2 \left|1-\delta_{kji}\mu_{kji}\right|= 2 \left|\delta_{ijk}\mu_{ijk}\right|
\le c_2 \exp \left(-\frac{A}{c_1}\right),
\]
where 
\[
c_2=2 \max_{i,j,k}|\delta_{ijk}|.
\]
Note that here we used the inequality $|\log z|<2|1-z|$, valid for $z\in\C$ if the right hand side is
$< 0.795$. This is satisfied, because even the reduced bound for $A$ will be about 45.
Also, in all our examples
$\log\left|\frac{\xi^{(k)}-\xi^{(j)}}{\xi^{(k)}-\xi^{(i)}}\right|$ 
was linearly dependent on 
$\log\left|\frac{\ve_1^{(j,k)}}{\ve_1^{(i,k)}}\right|,
\ldots,\log\left|\frac{\ve_5^{(j,k)}}{\ve_5^{(i,k)}}\right|$,
therefore we have to deal with an inequality of the form 
\begin{equation}
\Lambda=
\left|a_1\log\left|\frac{\ve_1^{(j,k)}}{\ve_1^{(i,k)}}\right|+\ldots 
+a_5\log\left|\frac{\ve_5^{(j,k)}}{\ve_5^{(i,k)}}\right|
\right|\le c_2 \exp \left(-\frac{A}{c_1}\right),
\label{linf}
\end{equation}
(with an obvious transformation of $a_1,\ldots,a_5$, depending on the example).

Applying a Baker-type lower estimate for $\Lambda$ 
(for example the estimates of \cite{bw}) we obtain
\[
\lambda<\exp(-C \log A),
\]
with a large constant $C$. Comparing it with the upper estimate (\ref{linf})
we conclude with an upper bound $A_B$ for $A$. In our example this was about $10^{34}$.

\section{Reduction}

Let $H$ be a constant. Consider the lattice spanned by the columns of the 6 by 5 matrix
\[
\left(
\begin{array}{ccccc}
1& 0&0&0&0 \\
0&1&0&0&0\\
\vdots&\vdots&\vdots&\vdots&\vdots\\
0&0&0&0&1\\
H\log\left|\frac{\ve_1^{(j,k)}}{\ve_1^{(i,k)}}\right|&
H\log\left|\frac{\ve_2^{(j,k)}}{\ve_2^{(i,k)}}\right|&
H\log\left|\frac{\ve_3^{(j,k)}}{\ve_3^{(i,k)}}\right|&
H\log\left|\frac{\ve_4^{(j,k)}}{\ve_4^{(i,k)}}\right|&
H\log\left|\frac{\ve_5^{(j,k)}}{\ve_5^{(i,k)}}\right|
\end{array}
\right)
\]
We apply a special case of Lemma 2.4 of \cite{book}:
\begin{lem}
If $a_1,\ldots,a_5$ satisfy (\ref{linf}), $A=\max_{1\le h\le 5}|a_h|\le A_0$
and $H$ is large enough so that the first vector $b_1$ of the LLL-reduced basis 
of the above lattice satisfies
\[
|b_1|\ge \sqrt{96}\cdot A_0,
\]
then
\[
A\le c_1(\log H+\log c_2-\log A_0).
\]
\end{lem}

The original bound obtained in the previous section is reduced in several subsequent steps.
A typical sequence is the following, starting with $A_0=10^{34}$:
\[
\begin{array}{|c|c|c|}
\hline
A_0&H& \text{reduced bound}\\\hline
10^{34}&10^{174}&344\\\hline
344&10^{17}&59\\\hline
59&10^{13}&45\\\hline
45&2\cdot 10^{12}&43\\\hline
\end{array}
\]
Using such large values, we must be very cautious with the accuracy of our calculations.
We used 250 digits. In the first steps the reduction is very efficient,
after a few steps the reduced bound is not smaller than the original bound,
then we stop the procedure.

\section{Sieving}

In our case we only have 5 exponents instead of 
9 and 14, respectively, in  \cite{quintic} and \cite{sextic}. 
However, if the reduced bound is $A_R$, then the total number of possible
exponent vectors $(a_1,\ldots,a_5)$ to check is $(2A_R+1)^5$, which, in case $A_R=45$, is more than
$6\cdot 10^9$. In this section we introduce a new idea, which helps to overcome this
difficulty, and can be applied in several similar cases.

The first idea is to use sieving, which turned out to be very efficient 
in a previous work \cite{sieve}. 

The main idea is to utilize the fact, that in all our examples one of the
fundamental units comes from the quadratic subfield.

Consider an equation of the form (\ref{xxijk}). As mentioned above, in all of our examples
one of the fundamental units was from the quadratic subfield $\Q(\sqrt{5})$.
Assume this is $\ve_1$.
In the terms of type
\[
\frac{(\xi^{(i)}-\xi^{(j)})\ell_{i,j}(x)}{(\xi^{(i)}-\xi^{(k)})\ell_{i,k}(x)} =
\frac{\xi^{(i)}-\xi^{(j)}}{\xi^{(i)}-\xi^{(k)}}
\prod_{1\le h\le 5}\left(\frac{\ve_h^{(i,j)}}{\ve_h^{(i,k)}}\right)^{a_h}
\]
it may happen, that $\ve_1^{(i,j)}=\ve_1^{(i,k)}$.
In Section \ref{str} we explicitly selected those conjugates $(i,j)$ 
which leave the elements of the quadratic subfield unchanged.
The problem is, that there exists no triple $(i,j,k)$, such that 
all of $(i,j),(j,k),(k,i)$ leaves the quadratic subfield fixed.
Therefore identities of type (\ref{xxijk}) can not be used to get rid of $\ve_1$.

But it turns out from (\ref{f2}), that all of 
\[
(1,2),(2,3),(3,4),(4,5),(5,1)
\]
has this property. We have 
\[
(\xi^{(1)}-\xi^{(2)})\ell_{1,2}(x)
+(\xi^{(2)}-\xi^{(3)})\ell_{2,3}(x)
+(\xi^{(3)}-\xi^{(4)})\ell_{3,4}(x)
\]
\[
+(\xi^{(4)}-\xi^{(5)})\ell_{4,5}(x)
+(\xi^{(5)}-\xi^{(1)})\ell_{5,1}(x)=0
\]
whence
\[
\frac{(\xi^{(1)}-\xi^{(2)})\ell_{1,2}(x)}{(\xi^{(1)}-\xi^{(5)})\ell_{1,5}(x)}
+
\frac{(\xi^{(2)}-\xi^{(3)})\ell_{2,3}(x)}{(\xi^{(1)}-\xi^{(5)})\ell_{1,5}(x)}
\]
\begin{equation}
+
\frac{(\xi^{(3)}-\xi^{(4)})\ell_{3,4}(x)}{(\xi^{(1)}-\xi^{(5)})\ell_{1,5}(x)}
+
\frac{(\xi^{(4)}-\xi^{(5)})\ell_{4.5}(x)}{(\xi^{(1)}-\xi^{(5)})\ell_{1,5}(x)}
=1.
\label{4term}
\end{equation}
In all fraction of type
\[
\frac{(\xi^{(i)}-\xi^{(j)})\ell_{i,j}(x)}{(\xi^{(1)}-\xi^{(5)})\ell_{1,5}(x)}
=
\frac{\xi^{(i)}-\xi^{(j)}}{\xi^{(1)}-\xi^{(5)}}
\left(\frac{\ve_1^{(i,j)}}{\ve_1^{(1,5)}}\right)^{a_1}
\left(\frac{\ve_2^{(i,j)}}{\ve_2^{(1,5)}}\right)^{a_2}
\left(\frac{\ve_3^{(i,j)}}{\ve_3^{(1,5)}}\right)^{a_3}
\left(\frac{\ve_4^{(i,j)}}{\ve_4^{(1,5)}}\right)^{a_4}
\left(\frac{\ve_5^{(i,j)}}{\ve_5^{(1,5)}}\right)^{a_5}
\]
the $\frac{\ve_1^{(i,j)}}{\ve_1^{(1,5)}}=1$, whence
for $(i,j)=(1,2),(2,3),(3,4),(4,5)$ we have
\[
\frac{(\xi^{(i)}-\xi^{(j)})\ell_{i,j}(x)}{(\xi^{(1)}-\xi^{(5)})\ell_{1,5}(x)}
=
\frac{\xi^{(i)}-\xi^{(j)}}{\xi^{(1)}-\xi^{(5)}}
\left(\frac{\ve_2^{(i,j)}}{\ve_2^{(1,5)}}\right)^{a_2}
\left(\frac{\ve_3^{(i,j)}}{\ve_3^{(1,5)}}\right)^{a_3}
\left(\frac{\ve_4^{(i,j)}}{\ve_4^{(1,5)}}\right)^{a_4}
\left(\frac{\ve_5^{(i,j)}}{\ve_5^{(1,5)}}\right)^{a_5}.
\]
This yields, that testing equations of the form (\ref{4term}),
instead of the usual identity (\ref{xxijk}),
we have to calculate 4 fractions for $(2A_R+1)^4$ exponent tuples,
instead of 2 fractions for $(2A_R+1)^5$ exponent tuples.
In case $A_R=45$, it yields calculating only 274299844
fractions, instead of 12480642902
($2.7\cdot 10^8$ instead of $1.2\cdot 10^{10}$).
Finally, it is important, that the integer arithmetic of Maple is very fast,
calculating $2.7\cdot 10^8$ fractions can be performed within a few minutes,
which makes our method feasible, in contrast with \cite{quintic} and \cite{sextic}.
Also, it is much simpler to implement this procedure.

We choose primes $p$, such that the polynomial 
$f(x)=x^{5}-m$  splits into linear factors modulo $p$, that is
\[
f(x)\equiv \prod_{1\le i\le 5} (x-m_{i})\;\; (\bmod\; p),
\]
with some $m_{i}\in\Z$. 
This yields $\xi^{(i)}\equiv m_{i}$ modulo a prime ideal $\frak{p}$
of the normal closure of $K=\Q(\xi)$, lying above $p$.
Then $\xi^{(1)}/\xi^{(2)}\equiv m_1/m_1 \;(\bmod \; \frak{p})$ gives a primitive fifth root of unity,
which allows to determine $\xi^{(j)}\;(\bmod \; \frak{p})$ $(2\le j\le 5)$ from 
$\xi^{(1)}$, similarly, as in Section \ref{str}. Then we obtain 
$\xi^{(i)}+\xi^{(j)},\xi^{(i)}-\xi^{(j)}$, and the corresponding conjugates of the units
modulo $\frak{p}$. Assume that 
\[
\frac{\xi^{(i)}-\xi^{(j)}}{\xi^{(1)}-\xi^{(5)}}\equiv s_{ij15}\;(\bmod \;\frak{p}),\;\;\;
\frac{\ve_h^{(i,j)}}{\ve_h^{(1,5)}}\equiv e_{hij15}\;(\bmod \;\frak{p}),
\]
Then (\ref{4term}) implies
\[
s_{1215}\prod_{h=2}^4 e_{h1215}^{a_h}
+s_{2315}\prod_{h=2}^4 e_{h2315}^{a_h}
+s_{3415}\prod_{h=2}^4 e_{h3415}^{a_h}
+s_{4515}\prod_{h=2}^4 e_{h4515}^{a_h}
\equiv 1\;(\bmod \;\frak{p}).
\]
But these are all integers in $\Z$, hence the congruence is valid also modulo the prime $p$:
\[
s_{1215}\prod_{h=2}^4 e_{h1215}^{a_h}
+s_{2315}\prod_{h=2}^4 e_{h2315}^{a_h}
+s_{3415}\prod_{h=2}^4 e_{h3415}^{a_h}
+s_{4515}\prod_{h=2}^4 e_{h4515}^{a_h}
\equiv 1\;(\bmod \;p).
\]
This is the modulo $p$ congruence that we test for $-A_R\le a_2,a_3,a_4,a_5\le A_R$.
If a tuple $(a_2,a_3,a_4,a_5)$ survives the test, then we test it modulo another suitable prime.
All together we used five primes, and only very few exponent tuples survived all of them.

If $(a_2,a_3,a_4,a_5)$ passed these tests, then we substituted $(a_2,a_3,a_4,a_5)$ into
an identity of the form (\ref{xxijk}), involving also $a_1$. We checked the identity
for $-A_R\le a_1 \le A_R$ modulo some primes $p$, similarly as above.

If $(a_1,a_2,a_3,a_4,a_5)$ passed, we explicitly calculated 
$\nu_{ij}=\prod_{h=1}^5 (\ve_h^{(ij)})^{a_h}$ (cf. (\ref{spd}))
for four pairs $i,j$, such that the system of linear equations
\[
\ell_{ij}(x)=\nu_{ij}
\]
is uniquely solvable in $x_1,x_2,x_3,x_4$. If the solutions are integers, 
we checked, if the index of $\alpha$ in (\ref{ax}) is indeed equal to 1.

\section{Computational results}

Using our algorithm, we calculated all inequivalent generators 
of power integral bases in ten pure quintic fields
generated by a root of an irreducible monogenic polynomial $f(x)=x^5-m$ with smallest positive $m$.

In the following tables we display $m$, the bound $A_B$ obtained by Baker's method,
the bound $A_R$ obtained after the reduction.
We include the primes that we used for sieving, 
the total number of exponent tuples to be tested
and the number of exponent tuples, that 
survived the modulo $p$ tests and were further tested by the next prime.
Finally, all $(x_1,x_2,x_3,x_4)$ are listed, such that 
$\alpha=x_0+x_1\xi+x_2\xi^2+x_3\xi^3+x_4\xi^4$ generates a power integral basis in $K$
($x_0\in\Z$ arbitrary).

\[
\begin{array}{l}
\hline
m=2, \;\; A_B=10^{34},\;\; A_R=54,\;\; \text{CPU time: 9 min}\\ \hline
\text{primes}: 151,241,251,431,571\\
\text{tested}: 141158161 \\
\text{survived}: 1869184, 15734, 254, 139, 136\\
\text{solutions}: (1, 0, 0, 0), (1, 1, 1, 1)\\
\hline
\end{array}
\]

\[
\begin{array}{l}
\hline
m=3, \;\; A_B=10^{34},\;\; A_R=47,\;\; \text{CPU time: 6 min}\\ \hline
\text{primes}: 41,431,491,661,761\\
\text{tested}:  81450625\\
\text{survived}: 3963567, 18458, 131, 44, 44\\
\text{solutions}: (1, 0, 0, 0)\\
\hline
\end{array}
\]

\[
\begin{array}{l}
\hline
m=5, \;\; A_B=10^{34},\;\; A_R=47,\;\; \text{CPU time: 5 min}\\ \hline
\text{primes}: 31,191,251,271,601\\
\text{tested}: 81450625 \\
\text{survived}: 5254688, 55052, 454, 9, 6\\
\text{solutions}: (1, 0, 0, 0)\\
\hline
\end{array}
\]

\[
\begin{array}{l}
\hline
m=6, \;\; A_B=10^{34},\;\; A_R=41,\;\; \text{CPU time: 3 min}\\ \hline
\text{primes}: 31,101,191,281,421\\
\text{tested}:  47458321\\
\text{survived}: 2936318, 58581, 633, 13, 8\\
\text{solutions}: (1, 0, 0, 0), (1, -1, -1, 0)\\
\hline
\end{array}
\]

\[
\begin{array}{l}
\hline
m=10, \;\; A_B=10^{34},\;\; A_R=43,\;\; \text{CPU time: 4 min}\\ \hline
\text{primes}: 11,101,251,281,521\\
\text{tested}: 57289761 \\
\text{survived}:  9894424, 195800, 1538, 15, 3\\
\text{solutions}: (1, 0, 0, 0)\\
\hline
\end{array}
\]

\[
\begin{array}{l}
\hline
m=11, \;\; A_B=10^{34},\;\; A_R=46,\;\; \text{CPU time: 5 min}\\ \hline
\text{primes}: 61,191,241,311,541\\
\text{tested}: 74805201 \\
\text{survived}: 2463570, 25691, 265, 48, 46\\
\text{solutions}: (1, 0, 0, 0), (1, -1, 0, 0)\\
\hline
\end{array}
\]

\[
\begin{array}{l}
\hline
m=13, \;\; A_B=10^{34},\;\; A_R=45,\;\; \text{CPU time: 5 min}\\ \hline
\text{primes}: 61,271,311,331,461\\
\text{tested}:  68574961\\
\text{survived}: 2263594, 16640, 127, 11, 10\\
\text{solutions}: (1, 0, 0, 0)\\
\hline
\end{array}
\]

\[
\begin{array}{l}
\hline
m=14, \;\; A_B=10^{34},\;\; A_R=46,\;\; \text{CPU time: 5 min}\\ \hline
\text{primes}: 41,61,101,191,211\\
\text{tested}:  74805201\\
\text{survived}: 3665415, 119766, 2389, 28, 1\\
\text{solutions}: (1, 0, 0, 0)\\
\hline
\end{array}
\]

\[
\begin{array}{l}
\hline
m=15, \;\; A_B=10^{34},\;\; A_R=43,\;\; \text{CPU time: 4 min}\\ \hline
\text{primes}: 211,241,311,541,761\\
\text{tested}:  57289761\\
\text{survived}: 543159, 4619, 17, 1, 1\\
\text{solutions}: (1, 0, 0, 0)\\
\hline
\end{array}
\]

\[
\begin{array}{l}
\hline
m=17, \;\; A_B=10^{34},\;\; A_R=46,\;\; \text{CPU time: 5 min}\\ \hline
\text{primes}: 101,181,491,601,701\\
\text{tested}:  74805201\\
\text{survived}: 1477946, 16109, 80, 1, 1\\
\text{solutions}: (1, 0, 0, 0)\\
\hline
\end{array}
\]

\section{Remarks on the calculations}

The units of the degree 10 fields were calculated by Magma \cite{magma}, all other calculations
were performed by Maple \cite{maple}. The codes were running on a usual PC with 2.5 GHz processor.
For the reduction procedure we used 250 digits accuracy.
In addition to the modulo $p$ test of the possible exponent tuples, all
remaining calculations took a negligible amount of CPU time, a few seconds.

\section{Conclusion}
\label{conc}
Introducing a new idea for testing small exponent tuples, we succeeded 
to present an algorithm to determine all generators of power integral bases in pure quintic fields.

This idea will certainly have several further applications.

For simplicity we presented our method for monogenic polynomials $x^5-m$. 
However, our procedure can obviously be extended to the case of non-monogenic polynomials
$f(x)=x^5-m$. In that case instead of (\ref{ax}) we represent $\alpha\in\Z_K$
in the form 
\[
\alpha=\frac{x_0+x_1\xi+x_2\xi^2+x_3\xi^3+x_4\xi^4}{d}
\]
with a common denominator $d\in\Z$. This implies that 
\[
\prod_{1\le i< j\le 5}|\ell_{i,j}(x)|=\frac{d^{10}}{I(\xi)}=N,
\]
that is, instead of being a unit in $L_{i,j}$, the $\ell_{i,j}(x)$
is an algebraic integer of norm $N$ in $L_{i,j}$.
Then we have to determine all non-associated integers $\gamma^{(i,j)}$ of $L_{i,j}$ of norm $N$,
\[
\ell_{i,j}(x)=\pm \gamma^{(i,j)} \prod_{h=1}^5 (\ve_h^{(i,j)})^{a_h},
\]
and the above procedure must be performed for all possible $\gamma$.

\end{document}